\documentclass[12pt]{article}
\usepackage{amsfonts}
\usepackage{graphicx}
\usepackage[latin1]{inputenc}
\usepackage[english]{babel}
\usepackage{epsfig}
\usepackage[dvips]{color}
\input xy \xyoption{all}
\usepackage{amssymb,latexsym}
\usepackage{euscript}

\usepackage{amsmath}
\usepackage{amscd}

\newtheorem{teo}{\bf Theorem}[section]
\newtheorem{prop}[teo]{\bf Proposition}
\newtheorem{ejem}[teo]{\bf Example}
\newtheorem{lema}[teo]{\bf Lemma}
\newtheorem{obser}[teo]{\bf Remark}
\newtheorem{defi}[teo]{\bf Definition}
\newtheorem{algo}[teo]{\bf Algorithm}
\newtheorem{coro}[teo]{\bf Corollary}

\newenvironment{Demo}{Proof.}

\def\Hom{\mathop{\rm Hom}\nolimits}
\def\rad{\mathop{\rm rad}\nolimits}

\def\gldim{\mathop{\rm gldim}\nolimits}

\def\qed{\hfill \mbox{$\square$}}

\def\Hom{\mathop{\rm Hom}\nolimits}
\def\Ext{\mathop{\rm Ext}\nolimits}

\def\soc{\mathop{\rm soc}\nolimits}
\def\pd{\mathop{\rm pd}\nolimits}
\def\id{\mathop{\rm id}\nolimits}

\def\dim{\mathop{\rm dim_k}\nolimits}

\def\top{\mathop{\rm top}\nolimits}

\def\h0{{\Hom_{A_0}}}

\begin{document}

\title{\large \bf The Higher Relation Bimodule}
\author{Ibrahim Assem, M. Andrea Gatica, Ralf Schiffler}
\date{}

\maketitle

\begin{abstract}

Given a finite dimensional algebra $A$ of finite global dimension, we consider the trivial extension of $A$ by the $A-A$-bimodule $\oplus_{i\ge 2} \Ext^2_A(DA,A)$, which we call the higher relation bimodule. We first give a recipe allowing to construct the quiver of this trivial extension in case $A$ is a string algebra and then apply it to prove that, if $A$ is gentle, then the tensor algebra of the higher relation bimodule is gentle.

\end{abstract}

\section{Introduction}

The objective of this paper is to describe a new class of algebras, which we call higher relation extensions. Our motivation comes from the study of cluster-tilted algebras, introduced by Buan, Marsh and Reiten in \cite{BMR}, and in \cite{CCS} for type $\mathbb{A}$. Indeed, it was shown in \cite{ABS} that an algebra $A$ is cluster-tilted if and only if there exists a tilted algebra $C$ such that $A$ is isomorphic to the trivial extension of $C$ by the $C-C$-bimodule $\Ext^2_C(DC,C)$. Moreover, a recipe for constructing the quiver of this trivial extension was given in \cite[Theorem 2.6]{ABS}. The proof of the latter result rests on the fact that tilted algebras have global dimension 2.

Here, we consider the more general case of an algebra $A$ having an arbitrary finite global dimension and consider its trivial extension by the bimodule $\bigoplus_{i\ge 2} \Ext^i_A(DA,A)$, which we call the higher relation bimodule. We believe that this class of algebras, which we call higher relation extensions, will be useful in the study of $m$-cluster-tilted algebras (see \cite{FPT} \cite{B}). Our first objective is to describe the ordinary quiver of the higher relation extension of $A$ in the case where $A$ is a string algebra in the sense of Butler and Ringel \cite{BR}. We also assume that the quiver of $A$ is a tree. This is no restriction, because the universal cover of a string algebra is a string tree \cite{G}. Our theorem reads as follows.

\begin{teo}\label{thmA}
Let $A=kQ/I$ be a string tree algebra. Then there exist two sequences $(c_\ell),(z_\ell)$ of points of $Q$ such that the arrows in the quiver of the higher relation extension are exactly those of $Q$ plus one additional arrow from each $z_\ell $ to $c_\ell$.
 \end{teo}

Our proof is constructive, in the sense that we give an algorithm allowing to construct explicitly the sequences $(c_\ell )$ and $(z_\ell) $ and thus the quiver of the higher relation extension.

We then consider the particular case where $A$ is a gentle algebra. Gentle algebras form an important subclass of the class of string algebras. Part of their importance comes from the fact that this subclass is stable under derived equivalences \cite{SZ}. While, as we show, the higher relation extension algebra of a gentle algebra is monomial but not necessarily gentle, we prove using our Theorem \ref{thmA}, that the tensor algebra of the higher relation bimodule is gentle.

\begin{teo}
\label{thmB}
Let $A=kQ/I$ be a gentle algebra, then the tensor algebra of the higher relation bimodule $\bigoplus_{i\ge 2} \Ext^i_A(DA,A)$ is gentle.
\end{teo}

The paper is organised as follows. In section 2, we fix the notation and recall some facts and results about string and gentle algebras. Section 3 is devoted to the computation of projective resolutions and injective coresolutions of uniserial modules over a string algebra. We study the top of the higher extension bimodule in section 4 and we prove Theorem \ref{thmA} in section 5. Sections 6 and 7 are devoted to the case of gentle algebras.

\section{ Preliminaries}

\subsection{ Notation}

Throughout this paper, algebras are basic and connected finite dimensional algebras over an algebraically closed field $k$. Given an algebra $A$, there always exists a (unique) quiver $Q=(Q_0,Q_1)$ and (at least) an isomorphism $A\cong kQ/I$, where $kQ$ is the path algebra of $Q$, and $I$ is an admissible ideal of $kQ$, see, for instance, \cite{ASS}. Such an isomorphism is called a {\bf presentation} of the algebra. Given an algebra $A$, we denote by $\textup{mod}\,A$ the category of finitely generated right $A$-modules, and by $D=\Hom_k(-,k)$ the standard duality between $\textup{mod}\,A$ and $\textup{mod}\,A^{op}$. For a point $x $ in the quiver $Q$ of $A$, we denote by $P(x),I(x),S(x),e_x$ respectively, the corresponding indecomposable projective module, injective module, simple module and primitive idempotent. We recall that a module $M$ can be equivalently considered as a bound quiver representation $M=(M_i,M_\alpha)_{i\in Q_0,\alpha\in Q_1}$. The projective, or injective, dimension of a module $M$ is denoted by $\pd M$, or $\id M$, respectively. The global dimension of $A$ is denoted by $\gldim A$. For facts about the category $\textup{mod}\, A$, we refer the reader to \cite{ARS} or \cite{ASS}.

\subsection {\bf Trivial extensions}

Let $A$ be an algebra and $M$ an $A-A$-bimodule . The 
trivial extension of $A$ by $M$ is the algebra $ A \ltimes M$ with
underlying $k$-vector space
$$A \oplus M = \{(a,m) | a \in A, \, m \in M \} $$ and the multiplication defined by
$$(a,m) . (a',m')= (a.a', am'+ma')$$
for $a,a' \in A$ and $m,m' \in M$.

For instance, an algebra $A$ is cluster-tilted if and only if there exists a tilted algebra $C$ such that $A$ is the trivial extension of $C$ by the so-called relation bimodule $\Ext^2_C(DC,C)$, see \cite{ABS}.

The ordinary quiver of a trivial extension is computed as follows (see, for instance, \cite{ABS}): let $M$ be an $A-A$ bimodule, then
the quiver $Q_{A \ltimes M}$ of $A \ltimes M$ is
given by

\begin{itemize}
\item[1)] $ \bigr( Q_{A \ltimes M}\bigl)_0= (Q_A)_0$
\item[2)] For $z,c \in (Q_A)_0$, the set of arrows in $Q_{A \ltimes M}$ from $z$
to $c$  equals the set of arrows in $Q_A$ from $z$ to $c$ plus

$$\dim  \displaystyle \frac{e_z M e_c}{e_z M (\rad A) e_c + e_z (\rad A) M e_c }$$
  additional arrows from $z$ to
$c$.

The latter arrows are called {\bf new} arrows, while the former are the {\bf old} arrows.
\end{itemize}

\subsection {String algebras}

Recall from \cite{BR} (see also \cite{WW}) that an algebra $A$ is called a {\bf string algebra} if
there exists a presentation $A=kQ/I$ 
(called {\bf a string presentation}) such that:
\begin{itemize}
\item[S1)] $I$ is generated by a set of paths (thus $A$ is monomial).
\item [S2)] Each point in $Q$ is the source of at most two
arrows and the target of at most two arrows.
\item[ S3)] For an arrow $\alpha$, there is at most one arrow
$\beta$ and at most one arrow $\gamma$ such that $\alpha \beta
\notin I$ and $\gamma \alpha \notin I$.
\end{itemize}

Whenever we deal with a string algebra $A$, we always assume that it is given by a string presentation $A=kQ/I$. We assume moreover that the relations (that is, the generators of $I$) are of minimal length.

A reduced walk $\omega$ in $Q$ is
called a {\bf string} if it contains no zero relations.To each string $\omega$ in $Q$, we can associate a
 so-called {\bf string module} \cite{BR} in the following way. If
$\omega $ is the stationary path at $j$, then $M(\omega)=S(j)$. Let $\omega = \omega_1 \omega_2 \cdots
\omega_t$ be a string, with each $\omega_i$ an arrow or the
inverse of an arrow. For each $i$ such that $0 \leq i \leq t$, let $V_i=k$; and for $1 \leq i \leq t$, let $V_{\omega_i}$ be the identity map sending $x
\in V_i$ to $x \in V_{i+1}$ if $\omega_i$ is an arrow and otherwise
the identity map sending $x \in V_{i+1}$ to $x \in V_i$. The string
module $M(\omega)$ is then defined as follows: for each $j \in Q_0$,
$M(\omega)_j$ is the direct sum of the vector spaces $V_i$ such that
the source of  $\omega_i$ is $j$ if $j$ appears in $\omega$, and otherwise
$M(\omega)_j=0$; for each $\alpha \in Q_1$, $M(\omega)_{\alpha}$ is the
direct sum of the maps $V_{\omega_i}$ such that $\omega_i=\alpha$ or
$\omega_i^{-1}= \alpha$ if $\alpha$ appears in $\omega$, and otherwise
$M(\omega)_{\alpha}=0$.

A non-zero path $\omega$ in $Q$ for $a$ to $b$ will sometimes be denoted by 
$ [a,b] $, whenever there is no ambiguity. Then, the corresponding string module is denoted by $M(\omega)=M([a,b])$.

We also recall that the endomorphism ring of a projective module over a string tree algebra $A$ (a full subcategory of $A$) is also a string tree algebra.

\subsection {Gentle algebras} 

Recall from \cite{AS} that 
a string algebra $A=kQ/I$ is called {\bf gentle} if in 
addition to $(S1),\, (S2), \, (S3)$, the bound quiver $(Q,I)$ satisfies:
\begin{itemize}
\item[G1)] For an arrow $\alpha$, there is at most one arrow
$\beta$ and at most one arrow $\gamma$ such that $\alpha \beta \in
I$ and $\gamma \alpha \in I$.

\item[G2)] $I$ is quadratic (that is, $I$ is generated by paths of length 2).
\end{itemize}

For instance, cluster-tilted algebras of types $\mathbb{A}$ and
$\tilde{ \mathbb{A}}$ are gentle \cite{ABCP}.

\section{ \large \bf Resolutions of uniserial modules}

In this section, we compute minimal projective resolutions of an injective module, and dually minimal injective coresolutions of a projective module over a string algebra. Throughout, we let $A=kQ/I$ be a string presentation.

\begin{defi} Let $[x_0,y_0]$ be a non-zero path from $x_0$ to $y_0$ in $Q$. We define inductively the right maximal sequence of $[x_0,y_0]$ as follows. This is a finite sequence of non-zero paths $[x_{i_1i_2 \cdots i_t}, y_{i_1i_2
\cdots i_t}]$ with $i_1=0$ and $i_j \in \{ 0,1\}$ such that 

\begin{itemize}

\item[1)] Let $[x_0,y_{00}]$,$[x_0,y_{01}]$ be the maximal non-zero paths
starting at $x_0$ (where we agree that $ [x_0,y_0]$ is contained in
$[x_0,y_{00}]$).

\noindent Then we set 
$$[x_{00},y_{00}]= [x_0,y_{00}] \backslash [x_0,y_0]$$
and
$$[x_{01},y_{01}]= [x_0,y_{01}] \backslash [x_0,y_0] =
[x_{0},y_{01}] \backslash \{x_0\}.$$

\item[2)] Inductively, assume that $[x_{0i_2 \cdots i_{t-1}}, y_{0i_2
\cdots i_{t-1}}]$ has been defined. Let $[x_{0i_2 \cdots i_{t-1}}, y_{0i_2 \cdots i_{t-1}0}]$
and $[x_{0i_2 \cdots i_{t-1}}, y_{0i_2 \cdots i_{t-1}1}]$ be the
maximal non-zero paths starting at $x_{0i_2 \cdots i_{t-1}}$, where
we agree that $[x_{0i_2 \cdots i_{t-1}}, y_{0i_2 \cdots i_{t-1}}] $
is contained in $[x_{0i_2 \cdots i_{t-1}}, y_{0i_2 \cdots
i_{t-1}0}]$.

\noindent Then we set 
$$[x_{0i_2 \cdots i_{t-1}0}, y_{0i_2 \cdots i_{t-1}0}]= [x_{0i_2
\cdots i_{t-1}}, y_{0i_2 \cdots i_{t-1}0}] \backslash [x_{0i_2
\cdots i_{t-1}}, y_{0i_2 \cdots i_{t-1}}]$$
and
$$[x_{0i_2 \cdots i_{t-1}1}, y_{0i_2 \cdots i_{t-1}1}]= [x_{0i_2
\cdots i_{t-1}}, y_{0i_2 \cdots i_{t-1}1}] \backslash \{ x_{0i_2 \cdots
i_{t-1}} \}.$$

\end{itemize}

\end{defi}

The left maximal sequence of a non-zero path is defined dually. However, we do it explicity for the convenience of the reader.

\begin{defi} Let $[r_0,s_0]$ be a non-zero path from $r_0$ to $s_0$ in $Q$. We define inductively the left maximal sequence of $[r_0,s_0]$ as follows. This is a finite sequence of non-zero paths $[r_{i_1i_2 \cdots i_t}, s_{i_1i_2 \cdots i_t}]$ with $i_1=0$ and  $i_j \in \{ 0,1\}$ such that

\begin{itemize}
\item[1)] Let $[r_{00},s_0]$,$[r_{01},s_0]$ be the maximal non-zero paths
ending at $s_0$, where we agree that $ [r_0,s_0]$ is contained in $[r_{00},s_0]$. Then we set 
$$[r_{00},s_{00}]= [r_{00},y_0] \backslash [r_0,s_0]$$
and
$$[r_{01}, s_{01}]= [r_{01},y_0] \backslash [r_0,s_0] =
[r_{01},s_0] \backslash \{s_0\}.$$

\item[2)] Inductively, assume that $[r_{0i_2 \cdots i_{t-1}}, s_{0i_2
\cdots i_{t-1}}]$ has been defined. Let $[r_{0i_2 \cdots i_{t-1}}, s_{0i_2 \cdots i_{t-1}0}]$
and $[r_{0i_2 \cdots i_{t-1}}, s_{0i_2 \cdots i_{t-1}1}]$ be the
maximal non-zero paths ending at $s_{0i_2 \cdots i_{t-1}}$, where we
agree that $[r_{0i_2 \cdots i_{t-1}}, s_{0i_2 \cdots i_{t-1}}] $ is
contained in $[r_{0i_2 \cdots i_{t-1}}, s_{0i_2 \cdots i_{t-1}0}]$.
Then we set 
$$[r_{0i_2 \cdots i_{t-1}0}, s_{0i_2 \cdots i_{t-1}0}]= [r_{0i_2
\cdots i_{t-1}}, s_{0i_2 \cdots i_{t-1}0}] \backslash [r_{0i_2
\cdots i_{t-1}}, s_{0i_2 \cdots i_{t-1}}]$$ and
$$[r_{0i_2 \cdots i_{t-1}1}, s_{0i_2 \cdots i_{t-1}1}]= [r_{0i_2
\cdots i_{t-1}}, s_{0i_2 \cdots i_{t-1}1}] \backslash \{ s_{0i_2 \cdots
i_{t-1}} \}.$$

\end{itemize}
\end{defi}

\noindent Note that in both cases, some of the paths above might be empty and in this case, the points considered do not exist.

\medskip

Our first result follows directly from the above definitions.

\begin{teo} \label{uniserial} Let $A=kQ/I$ be a string algebra.

\begin{itemize}
\item[a)] If  $[x_0,y_0]$ is a non-zero path in $Q$ and 
$$ \cdots \rightarrow P_3 \rightarrow P_2 \rightarrow P_1
\rightarrow P_0 \rightarrow M[x_0,y_0] \rightarrow 0$$ is a minimal
projective resolution then,  for $l \geq 1$,

$$ P_{l-1}= \bigoplus P(x_{i_1i_2 \cdots i_l })$$

\noindent where the direct sum is taken over all $l$-tuples $(0,i_2, \cdots,i_l)$
such that $ i_k \in \{ 0,1\}$ for all $k$ with $2 \leq k \leq l$ and the
point $x_{0i_2 \cdots i_l }$ in definition 3.1 exists.

\item[b)] If $[r_0.s_0]$ is a non-zero path in $Q$ and 
$$ 0  \rightarrow M[r_0,s_0] \rightarrow I_0 \rightarrow I_1 \rightarrow I_2
\rightarrow I_3 \rightarrow \cdots $$ is a minimal injective
coresolution then,  for $l \geq 1$,

$$ I_{l-1}= \bigoplus I(s_{i_1i_2 \cdots i_l })$$

\noindent where the direct sum is taken over all $l$-tuples $(0,i_2, \cdots,i_l)$
such that $ i_k \in \{ 0,1\}$ for all $k$ with $2 \leq k \leq l$ and the
point $s_{0i_2 \cdots i_l }$ in definition 3.2 exists.

\end{itemize}
\end{teo}

\begin{Demo} We only prove a), since the proof of b) is dual.

\noindent Clearly, the projective cover of the uniserial module $M[x_0,y_0]$ is 
$P(x_0)$, whose support consists of 
the (at most two) maximal non-zero paths 
$[x_0,y_{01}]$ and $[x_0,y_{01}]$ starting at
$x_0$. Then,$$ \Omega^1 M[x_0.y_0]= M[x_{00},y_{00}] \oplus
M[x_{01},y_{01}]$$

\noindent where $x_{00}$ and $x_{01}$ are defined as above. The rest follows from an easy induction.  $\qed$

\end{Demo}

\bigskip

\begin{ejem} \label{example 1} Suppose that the string algebra is given by the
following bound quiver.

\[\xymatrix@R12pt@C12pt{ \scriptscriptstyle 1 \ar[r] \ar@/^2ex/@{.}[rr]& \scriptscriptstyle 3  \ar[dr] \ar[r] \ar@/^2ex/@{.}[rrr] \ar@/^1ex/@{.}[rdr]
\ar@/_4pc/@{.}[dddrrrrrru]&
 \scriptscriptstyle 4 \ar[r] &\scriptscriptstyle 5  \ar[r] & \scriptscriptstyle 6\\
\scriptscriptstyle 2 \ar[ur] \ar@/^2ex/@{.}[rr]& &
\scriptscriptstyle 7 \ar[r] \ar[dr]&  \scriptscriptstyle 8 & &
\scriptscriptstyle 11 \ar[r] \ar@/^2ex/@{.}[rr] &
 \scriptscriptstyle 12 \ar[r] & \scriptscriptstyle 13\\
& & &  \scriptscriptstyle 9
 \ar[r]\ar@/^2ex/@{.}[rru] & \scriptscriptstyle 10 \ar[r] \ar[ur]\ar@/^3ex/@{.}[rur] & \scriptscriptstyle 14 \ar[r] &
 \scriptscriptstyle 15 \ar[r] & \scriptscriptstyle 16 \ar[r] & \scriptscriptstyle 17 \\
 & & & & & & & & &}
\]

\medskip

 \noindent Here, and in the sequel, dotted lines indicate relations.

\noindent Considering the path $[x_0,y_0]=[3,9]$, the right maximal sequence is  
$$ [3,9]; [10,15], [4,5]; [16,17], [11,11]=\{11\}, [6,6]=\{6\};[12,12]=\{12\};[13,13]=\{13\}.$$

\noindent This sequence may be conveniently shown in the following diagram

\[\xymatrix@R12pt@C12pt{  [3,9] \ar@{-}[r]  \ar@{-}[dr]& [4,5]  \ar@{-}[r] & \{6\} & & &\\
& [10,15] \ar@{-}[r]  \ar@{-}[dr]& \{11 \}\ar@{-}[r] & \{12\} \ar@{-}[r] & \{13\} & \\  
 & &
[16,17] & & .&}
\]

\noindent The minimal projective resolution of $M[3,9]$ is the following (compare with the above diagram)

$$ \begin{array}{c}
0 \rightarrow  P(13) \rightarrow P(12) \rightarrow\\ \\ \rightarrow P(16) \oplus P(11)
\oplus P(6)   \rightarrow P(10) \oplus P(4) \rightarrow P(3)
\rightarrow M[3,9] \rightarrow 0,
\end{array}$$

\noindent where the morphisms are induced by the corresponding paths.

\noindent Similarly, taking $[r_0,s_0]=[3,9]$, the left maximal sequence is 

\[\xymatrix@R12pt@C12pt{   \{1 \} \ar@{-}[r]  & [3,9]&}
\]

\noindent from which we deduce the minimal injective coresolution
$$ 0 \rightarrow M[3,9] \rightarrow I(9) \rightarrow I(1) \rightarrow 0.$$

\end{ejem}

We are interested in computing resolutions of injective and projective indecomposable modules. These modules are usually not uniserial, neither are in general their first syzygy or cosyzygy, respectively. In order to apply Theorem   \ref{uniserial}, the next lemma shows that we must start from the second.

\begin{lema} 

\begin{itemize}
\item[a)]\label{syzygy} The second syzygy of an indecomposable
injective module is the direct sum of at most six uniserial modules.
\item[b)]\label{cosyzygy} The second cosyzygy of an indecomposable projective module is the direct sum of at most six uniserial modules.
\end{itemize}
\end{lema}

\begin{Demo}  Let $I(c)$ be an indecomposable
injective $A$-module.  If $I(c)$ is uniserial, then there is nothing to prove because of Theorem \ref{uniserial}. Otherwise, let $\top I(c)=S(a_0) \oplus S(a_1)$. Then the projective cover of $I(c)$ is $P(a_0) \oplus P(a_1)$. Let
$[a_i, b_i]$ and $[a_i, b'_i]$ be the two maximal non-zero paths
starting at $a_i$ (with $i=0,1$), where we agree that $[a_0,c]$ is contained in $[a_0,b'_0]$ and $[a_1,c]$ is contained in $[a_1,b'_1]$. Let $d_i$ be the
direct successor of $a_i$ on the path $[a_i,b_i]$ then
$$ \Omega^1 I(c)=M[d_0,b_0] \oplus M[d_1,b_1] \oplus M$$
\noindent where $M$ is an indecomposable module, usually non-uniserial, such that $\top M=S(c)$ and $\soc M= S(b'_0) \oplus S(b'_1)$.
Hence, the projective cover of $\Omega^1 I(c)$ is $P(d_0) \oplus
P(d_1) \oplus P(c)$, and $\Omega^2 I(c)$ is the direct sum of at
most six uniserial modules obtained as follows.

\noindent Let $[d_i,b_{i0}], [d_i,b_{i1}]$ be the maximal non-zero
paths starting in $d_i$ (with $i=0,1$), where we agree that $[d_i,b_i]
$ is contained in $[d_i,b_{i0}]$. Then let

$$[d_{i0},b_{i0}]= [d_i,b_{i0}]\backslash [d_i,b_i]$$

$$[d_{i1},b_{i1}]= [d_i,b_{i1}]\backslash [d_i,b_i]= [d_i,b_{i1}] \backslash \{d_i\}.$$

\noindent Let also $[c,c_0]$ and $[c,c_1]$ be the maximal non-zero
paths starting at $c$, where we agree that $[c,b'_0]$ is contained in $[c,c_0]$ and $[c,b'_1]$ is contained in $ [c,c_1].$ 

\noindent We let $$[c'_0,c_0]=[c,c_0]
\backslash [c,b'_0]$$ and $$[c'_1,c_1]=[c,c_1] \backslash [c,b'_1].$$

\noindent It is then clear that
\begin{eqnarray*} \Omega^2 I(c)& = &  M[d_{00},b_{00}] \oplus M[d_{01},b_{01}] \\
& \oplus &  M[d_{10},b_{10}] \oplus M[d_{11},b_{11}] \\
& \oplus &  M[c'_0,c_0] \oplus M[c'_1,c_1]
\end{eqnarray*}
which establishes $a)$. Statement $b)$ is dual. $\qed$

\end{Demo}

 \begin{coro}

\begin{itemize}
\item[a)] \label{MinProjResInjective} Let $I(c)$ be an indecomposable injective module  such that $\top (I(c))= S(a_0) \oplus S(a_1)$.
Then $I(c)$ has the following minimal projective resolution

\medskip

{$$\begin{array}{c} \dots \rightarrow \bigoplus_{j; (0, i_2,i_3)} P(x^j_{0 i_2 i_3}) \rightarrow \bigoplus_{j; (0,i_2)} P(x^j_{0 i_2})
\rightarrow \bigoplus_j P(x^j_{0}) \rightarrow\\ \\
 \rightarrow
 P(d_0)\oplus P(c)
\oplus P(d_1) \rightarrow P(a_0) \oplus P(a_1) \rightarrow
I(c)\rightarrow 0
\end{array}$$ }

\noindent with the morphisms induced by the paths, where  $\{ x_0^j \, | \, 1 \leq j \leq 6 \}= \{d_{00},d_{01},d_{10},d_{11}, c'_{0},c'_1 \}$ and $ i_j \in \{0,1 \}$.

\item[b)] Let $P(z)$ be an indecomposable projective module 
 such that $\soc (P(z))= S(w_0) \oplus S(w_1)$. Then $P(z)$ 
has the following minimal injective coresolution 
$$ \begin{array}{c} 0  \rightarrow P(z) \rightarrow   I(w_0) \oplus I(w_1) \rightarrow  I(v_1)\oplus I(z)
\oplus I(v_2) \rightarrow \bigoplus_j I(s^j_{0}) \rightarrow
\\ \\
\rightarrow \bigoplus_{j; (0,i_2)} I(s^j_{0 i_2}) \rightarrow  \bigoplus_{j; (0,i_2,i_3)} I(s^j_{0 i_2 i_3}) \rightarrow
 \dots\end{array} $$ 
\medskip

\noindent with the morphisms induced by the paths, where  $\{ s_0^j \, | \, 1 \leq j \leq 6 \}$ are as above and $ i_j \in \{0,1 \}$.

\end{itemize}
\end{coro}

\begin{Demo} This follows from Lemma \ref{syzygy} and Theorem \ref{uniserial}. $\qed$

\end{Demo}

\begin{coro} With the above notations 
\begin{itemize}
\item[a)] All the points $x^j_0,
x^j_{0i_2}, \cdots, x^j_{0i_2\dots i_l}$ are targets of
relations.

\item[b)] All the points $s^j_0,
s^j_{0i_2}, \cdots, s^j_{0i_2\dots i_l}$ are sources of
relations.
\end{itemize}
\end{coro}
\begin{Demo}
This follows from the construction of these
points. $\qed$ 
\end{Demo}

\section{ The top of the higher relation bimodule}

\begin{defi} Let $A$ be a finite dimensional algebra of finite
global dimension. The $A-A$-bimodule $ \bigl(
\bigoplus_{i \geq 2}\Ext^i_A(DA,A)\bigr) $ with the natural
action is called the {\bf higher relation bimodule}. 
 The trivial extension
$$A \ltimes \bigr( \bigoplus_{i \geq 2}\Ext^i_A(DA,A)\bigl)$$

\noindent of $A$ by its higher relation bimodule is called the {\bf higher relation extension} of $A$.

\end{defi}

If $\gldim A \leq 2$, then the higher relation extension of $A$
coincides with its relation extension, as defined in \cite{ABS}.

Our objective in this section is to construct the ordinary quiver of the higher relation extension of a string algera $A$ of finite global dimension.

As mentioned in the introduction, we also assume that the ordinary quiver $Q_A$ of $A$ is a tree. 

Let thus $A=kQ/I$ be a string algebra, with $Q$ a tree and $M$ be an $A-A$-bimodule.
We have $$\rad M=M(\rad A) + (\rad A)M $$ and then $$\top M= M/ [M(\rad A) + (\rad A)M].$$

If $M= \bigoplus_{i \geq 2}\Ext^i_A(DA,A)$, then, clearly, $\top M=  \bigoplus_{i \geq 2} \top \Ext^i_A(DA,A)$. In order to describe this top, we start by describing the modules $\top_A \Ext^i_A(I(c),A)$ and $\top \Ext^i_A(DA,P(z))_A$ for all points $c,z \in (Q_A)_0$.

In the following, we use the notation of section 3.

\begin{prop}\label {Ext(I,A)} Let $A=kQ/I$ be a string tree algebra and $l \geq 0$.
Then $\Ext_A^{l+2}(I(c), P(z)) \ne 0$ if and only if one of the following two conditions hold:
\begin{itemize}
\item[a)] there exists a non-zero path 
$\omega: z \rightsquigarrow x^j_{i_1i_2 \cdots i_{l+1}}$ not passing through
$x^j_{i_1i_2 \cdots i_l}$ and whose compositions with 
$x^j_{i_1i_2 \cdots i_{l+1}}\rightsquigarrow x^j_{i_1i_2 \cdots
i_{l+2}}$ are both zero.
\item[b)] $z=x^j_{i_1i_2 \cdots i_l}$ and $x^j_{i_1i_2 \cdots i_{l}0}, \, x^j_{i_1i_2 \cdots i_{l}1}$ both exist. In this case, a non-zero element is induced from the difference of the two paths $x^j_{i_1i_2 \cdots i_{l}}\rightsquigarrow x^j_{i_1i_2 \cdots
i_{l}0}$ and  $x^j_{i_1i_2 \cdots i_{l}}\rightsquigarrow x^j_{i_1i_2 \cdots
i_{l}1}$.
 
\end{itemize}
\end{prop}

\begin{obser}
Observe that in case $(b)$, we have the following situation

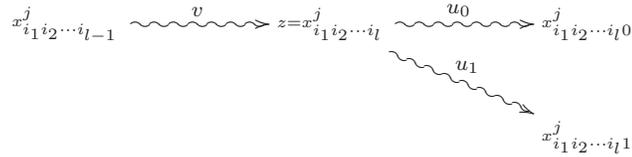
\[\xymatrix{ & \scriptscriptstyle x^j_{i_1i_2\cdots i_{l-1}}  \ar@{~>}^{v}[rr] 
& & \scriptscriptstyle z=x^j_{i_1i_2\cdots i_{l}} \ar@{~>}^{u_0}[rr] \ar@{~>}^{u_1}[drr] & &
 \scriptscriptstyle x^j_{i_1i_2\cdots i_{l}0} \\
& & & & & \scriptscriptstyle x^j_{i_1i_2\cdots i_{l}1} } \]

\noindent where $vu_0,vu_1$ are zero paths.
\end{obser}

\begin{Demo} 
Let 
$$ \cdots  \rightarrow \oplus P(x^j_{i_1 \cdots i_{l+2}}) \stackrel{d_{l+3}}{\to} \oplus P(x^j_{i_1 \cdots i_{l+1}}) \stackrel{d_{l+2}}{\to} \oplus P(x^j_{i_1 \cdots i_{l}})  \rightarrow \cdots \rightarrow P_{c} \rightarrow I(c) \rightarrow 0$$ 
be a minimal projective resolution of $I(c)$. Recall that the morphisms $d_k$ are induced from the paths in $Q$.  

\noindent If condition (a) holds then it follows from the definition of  $\Ext_A^{l+2}(I(c),
P(z))$ that $\omega$ induces a non-zero element in 
$\Ext_A^{l+2}(I(c),
P(z))$. 

\noindent If condition (b) holds, then $P_{l+2}=\oplus P(x^j_{i_1 \cdots i_{l+1}})$ has two indecomposable summands $P(x^j_{i_1 \cdots i_l0})$,  $P(x^j_{i_1 \cdots i_l1})$, whose images  $d(P(x^j_{i_1 \cdots i_l0}))$ and $d(P(x^j_{i_1 \cdots i_l1}))$ lie in the same indecomposable summand $P(x^j_{i_1 \cdots i_l})$ of $P_{l+1}$, together with two non-zero morphisms $\nu_i:P(x^j_{i_1 \cdots i_l i_{l+1}}) \rightarrow P(z)$ ($i_{l+1}=0,1$) such that there exist two morphisms $\gamma_i: P(x^j_{i_1 \cdots i_l}) \rightarrow P(z)$ with $\nu_i=\gamma_i d$.

\[\xymatrix{ P(x^j_{i_1 \cdots i_l0})\oplus P(x^j_{i_1 \cdots i_l1})\ar@<2pt>[r]^-{d}\ar@<-2pt>[d]_{\nu_1}\ar@<2pt>[d]^{\nu_2} 
& P(x^j_{i_1 \cdots i_l})\ar[dl]^{[\gamma_1,\gamma_2]} \\ P(z)
}\]

In this case $[\nu_1 \, \,- \nu_2]^t: P(x^j_{i_1 \cdots i_l0}) \oplus P(x^j_{i_1 \cdots i_l1}) \rightarrow P(z)$ does not factor through $P(x^j_{i_1 \cdots i_l})$ because $\dim \Hom (P(x^j_{i_1 \cdots i_l0}) \oplus P(x^j_{i_1 \cdots i_l1}),P(z))=2$ while $\dim \Hom (P(x^j_{i_1 \cdots i_l}),P(z))=1$, since the algebra is a tree algebra. This shows that $\Ext^{l+2}_A(I(c),P(z)) \ne 0$.

Conversely, suppose that 
$\Ext^{l+2}_A(I(c),P(z))$ contains a non-zero element $[f]$. Then $[f]$ is in the class of a morphism $f: \oplus P(x^j_{i_1 \cdots i_{l+1}}) \rightarrow P(z)$ such that $f d_{l+3}=0$. Since 
$A$ is a tree string algebra, there are at most two indecomposable summands on which $f$ is non-zero, because otherwise there are non-zero paths from $z$ to three points $x^j_{i_1 \cdots \i_{l+1}}$ and these induce a full subcategory of type $\mathbb{D}_4$ which contradicts the fact that $A$ is string. Thus we get a morphism $f: P(x^j_{i_1 \cdots i_{l+1}}) \oplus P(x^{j'}_{i'_1 \cdots i'_{l+1}}) \rightarrow P(z)$ which does not factor through $d_{l+2}$.
If $z=x^j_{i_1 \cdots i_l}$ then we must have $j=j'$, $i=i', \cdots, i_l=i'_{l}$ and $i_{l+1} \ne i'_{l+1}$.
Suppose $z \ne x^j_{i_1 \cdots i_l}$. 
If both non-zero paths $z \rightsquigarrow x^j_{i_1 \cdots i_{l+1}}$, $z \rightsquigarrow x^{j'}_{i'_1 \cdots i'_{l+1}}$ 
which induce $f$ pass through $x^j_{i_1 \cdots i_l}$
then we have a contradiction to $A$ being string. If one non-zero path $z \rightsquigarrow x^j_{i_1 \cdots i_{l+1}}$
passes through $x^j_{i_1 \cdots i_l}$
then the other satisfies condition $(a)$. Indeed, the composition with 
 $x^{j'}_{i'_1i'_2 \cdots i'_{l+1}}\rightsquigarrow x^{j'}_{i'_1i'_2 \cdots
i'_{l+2}}$ vanishes because our original path corresponds to an element of $\Ext^{l+2}_A(I(c),P(z))$.
Similarly, if $z \rightsquigarrow x^j_{i_1 \cdots i_{l+1}}$ does not pass through neither $x^j_{i_1 \cdots i_{l}}$ nor $x^{j'}_{i'_1 \cdots i'_{l+1}}$
then both paths satisfy condition $(a)$.
$\qed$
\end{Demo}

The following example ilustrates condition (b). 
\begin{ejem}
Let $A$ be given by the quiver

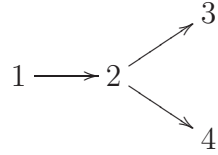
\[\xymatrix@R=10pt{
&&3\\
1\ar[r] &2\ar[ru]\ar[rd] \\
&&4
}\]
bound by $\textup{rad}^2A=0$. Then the minimal projective resolution of $I(1)$ is
\[0\to P(3)\oplus P(4) \to P(2) \to P(1)\to I(1)\to 0.
\]
Let $j_1:P(3)\to P(2) $ and $j_2:P(4)\to P(2)$ be the  canonical inclusions, then it is easily seen that the morphism 
\[ [j_1 \, \,\,\, \, -j_2 ]^t : P(3)\oplus P(4) \to P(2)
\]
induces a non-zero element of $\Ext^2_A(I(1),P(2)).$
\end{ejem}

\begin{coro}
Assume $A$ is a gentle tree algebra, then $\Ext_A^{l+2} (I(c),P(z))\ne 0$ if and only if  there exists a non-zero path 
$\omega: z \rightsquigarrow x^j_{i_1i_2 \cdots i_{l+1}}$ not passing through
$x^j_{i_1i_2 \cdots i_l}$ and whose compositions with 
$x^j_{i_1i_2 \cdots i_{l+1}}\rightsquigarrow x^j_{i_1i_2 \cdots
i_{l+2}}$ are both zero.
\end{coro}
\begin{Demo}
Indeed, if $A$ is gentle, then condition (b) cannot occur as shown in the remark preceding the proof. $\qed$
\end{Demo}

\begin{coro} \label {exists} \begin{itemize} 
\item[a)] Let $\omega: z \rightsquigarrow x^j_{i_1
i_2 \cdots i_{l+1}}$ be a non-zero path as in Proposition \ref{Ext(I,A)} a).
\begin{itemize}
\item[a1)] Assume that a point $x^j_{i_1 i_2 \cdots i_{l+2}}$ exists, then $w$ induces an element of the top of
$_A\Ext_A^{l+2}(I(c), A)$ if and only if $z$ is the starting point
of a relation of the form $\omega \omega '$ when $\omega': x^j_{i_1
i_2 \cdots i_{l+1}} \rightsquigarrow x^j_{i_1 i_2 \cdots i_{l+2}}$

\item[a2)] Assume that no point $x^j_{i_1 i_2 \cdots i_{l+2}}$
exists, then $w$ induces an element of the top of
$_A\Ext_A^{l+2}(I(c), A)$ if and only if $z= x^j_{i_1 i_2 \cdots
i_{l+1}} $ and $\omega $ is the stationary path in $z$.
\end{itemize}
\item[b)] In the situation of Proposition  \ref{Ext(I,A)} b), the class of the difference of the paths  $x_{i_1\ldots i_l}^j\rightsquigarrow x^j_{i_1 \ldots i_l0}$ and $x_{i_1\ldots i_l}^j\rightsquigarrow x^j_{i_1 \ldots i_l1}$ in $\Ext^{l+2}_A(I(c),P(z))$ lies in the top of  $_A\Ext_A^{l+2}(I(c), A)$ if and only if there are two minimal relations $z\rightsquigarrow x^j_{i_1\ldots i_l0}
\rightsquigarrow x^j_{i_1 \ldots i_l 0 i_{l+2}}
$ and $ z\rightsquigarrow x^j_{i_1 \ldots i_l 1}
\rightsquigarrow x^j_{i_1 \ldots i_l 1 i_{l+2}}$.
\end{itemize}\end{coro}

\begin{Demo} 
\begin{itemize}
\item[a1)] The morphism $f:P_{l+2} =\oplus P(x^j_{i_1\ldots i_li_{l+1}}) \rightarrow P(z)$ induced by
$\omega$ factors through $P(s)$ where $s$ is the source of a
relation ending at $x^j_{i_1 i_2 \cdots i_{l+2}}$ and such that $s$
lies on the path $\omega$. So, $f$ induces an element on the top of
$_A\Ext_A^{l+2}(I(c), A)$ if and only if $s=z$. 

\item[a2)] This follows from the fact that 
the morphism $f:P_{l+2}=\oplus P(x^j_{i_1\ldots i_li_{l+1}})
\rightarrow P(z)$ factors through the identity of $P(x^j_{i_1 i_2
\cdots i_{l+1}})$. 

\item[b)] Let $f$ be a representative of the class of the difference of paths  $x_{i_1\ldots i_l0}^j\rightsquigarrow x^j_{i_1 \ldots i_l}$ and $x_{i_1\ldots i_l 1}^j\rightsquigarrow x^j_{i_1 \ldots i_l}$ in $\Ext_A^{l+2}(I(c), P(x^j_{i_1 \ldots i_l}))$. Then $$f=[f_0 \, \, \,f_1 \, \,\,0]:P(x^j_{i_1\ldots i_l0})\oplus P(x^j_{i_1 \ldots i_l 1}) \oplus \overline{P}_{\ell+2}\longrightarrow P(x^j_{i_1 \ldots i_l}).$$ 
Suppose first that there is no relation  $z\rightsquigarrow x^j_{i_1\ldots i_l0}
\rightsquigarrow x^j_{i_1 \ldots i_l 0 i_{l+2}}
$.
Then any relation ending at $ x^j_{i_1 \ldots i_l 0 i_{l+2}}
$ must start at a successor $y$ of $z$. Therefore there exists $g:P( x^j_{i_1\ldots i_l0})\to P(y)$ such that $f_0$ factors through $g$, whence $ f=[hg \, \, \, f_1 \, \, \,0]$, for some morphism $h$. So $[f]$ is not in the top of $_A\Ext_A^{l+2}(I(c), A)$.

Conversely, if we have two minimal relations as in the statement, and $[f]$ is not in the top of $_A\Ext_A^{l+2}(I(c), A)$, then $[f]=[h][g]$ for some $[g]\in \Ext_A^{l+2}(I(c), P(y))$, which  is represented by a morphism \[g:P(x^j_{i_1\ldots i_l0})\oplus P(x^j_{i_1 \ldots i_l 1}) \longrightarrow P(y).\]
Then $y$ lies on the non-zero path $z\rightsquigarrow x_{i_1 \ldots i_l 0}^j$ or  $z\rightsquigarrow x_{i_1 \ldots i_l 1}^j$ (or both) and $y\ne z$.
But then $g\,d_{l+3,0}:P( x_{i_1 \ldots i_l 0 i_{l+2}}^j)\to P(y)$ is  non-zero, because it is given by the non-zero path $y\rightsquigarrow x_{i_1 \ldots i_l 0}^j \rightsquigarrow x_{i_1 \ldots i_l 0 i_{l+2}}^j$, and this contradicts the fact that $[g]$ belongs to $\Ext^{l+2}_A(I(c),P(y))$. $\qed$
\end{itemize}

\end{Demo}

We summarise the results in the theorem below.

{ For each point $c$ in a string algebra $A=kQ/I$, we compute the
minimal projective resolution of $I(c)$ given in Corollary
\ref{MinProjResInjective}. Then for all $l \geq 0$, the $l+2$-nd
term in the minimal projective resolution of $I(c)$ is given by
$P_{l+2}= \bigoplus_{j, (i_1,i_2, \cdots , i_{l+1}) } P(x^j_{i_1 i_2
\cdots i_{l+1}})$.  }

{Whenever the point $x^j_{i_1 i_2 \cdots i_{l+2}}$
exists, let $z^j_{i_1 i_2 \cdots i_{l+2}}$ be the source of
the relation ending in $x^j_{i_1 i_2 \cdots i_{l+2}}$ and passing
through $x^j_{i_1 i_2 \cdots i_{l+1}}$.}

For each $j \in \{ 1, \cdots, 6\}$ and for
each $l \geq 0$, define
$$\mathbf{Z}^j_{i_1 i_2 \cdots i_{l+1}} = \left\{\begin{array}{ll} z^j_{i_1 i_2 \cdots i_{l+1} 0} &\textup{if $x^j_{i_1 i_2 \cdots i_{l+1} 0}$  exists;} \\ 
x^j_{i_1 i_2 \cdots i_{l+1} } &\textup{otherwise,} \end{array}\right.$$
and let 
\[\zeta_{i_1\ldots i_{l+1}}^j = \left\{\begin{array}{ll}
[x^j_{i_1 \cdots i_l}, \mathbf{Z}^j_{i_1 i_2 \cdots i_{l+1}} ]
&\textup{if $x^j_{i_1i_2 \cdots i_{l}0}, \, x^j_{i_1i_2 \cdots i_{l}1}$
both exist;} \\
\left[x^j_{i_1 \cdots i_l}, \mathbf{Z}^j_{i_1 i_2 \cdots i_{l+1}}\right] \setminus \{ x^j_{i_1 \cdots i_l} \}
& \textup{otherwise.}
\end{array}\right.
\]
Dually, 
whenever the point $s^j_{i_1 i_2 \cdots i_{l+2}}$
exists, let $c^j_{i_1 i_2 \cdots i_{l+2}}$ be the target of
the relation starting in $s^j_{i_1 i_2 \cdots i_{l+2}}$ and passing
through $s^j_{i_1 i_2 \cdots i_{l+1}}$.
For each $j \in \{ 1, \cdots, 6\}$ and for
each $l \geq 0$, define
$$\mathbf{C}^j_{i_1 i_2 \cdots i_{l+1}} = \left\{\begin{array}{ll} c^j_{i_1 i_2 \cdots i_{l+1} 0} &\textup{if $s^j_{i_1 i_2 \cdots i_{l+1} 0}$  exists;} \\ 
s^j_{i_1 i_2 \cdots i_{l+1} } &\textup{otherwise,} \end{array}\right.$$
and let 

\[\Theta_{i_1\ldots i_{l+1}}^j = \left\{\begin{array}{ll}
[ \mathbf{C}^j_{i_1 i_2 \cdots i_{l+1}}, s^j_{i_1 \cdots i_l}]
&\textup{if $s^j_{i_1i_2 \cdots i_{l}0}, \, s^j_{i_1i_2 \cdots i_{l}1}$
both exist;}\\
\left[ \mathbf{C}^j_{i_1 i_2 \cdots i_{l+1}}, s^j_{i_1 \cdots i_l} \right] \setminus \{ s^j_{i_1 \cdots i_l} \}  
&\textup{otherwise.}
\end{array}\right.
\]

\begin{teo} \label{left top}Let $A=kQ/I$ be a string tree algebra and
 $l \geq 0$. The following are equivalent
\begin{itemize}
\item[{\rm (a)}] $ \Ext_A^{l+2}(I(c), P(z)) \ne 0$;
\item[{\rm (b)}]  there exists $j$ such that $z \in \zeta^j_{i_1 i_2 \cdots i_{l+1}}$;
\item[{\rm (c)}]  there exists $j$ such that $c \in \Theta^j_{i_1 i_2 \cdots i_{l+1}}$.
\end{itemize}
\end{teo}

\begin{Demo} The equivalence of (a) and (b) follows from Proposition \ref{MinProjResInjective} and from the definition of  $\zeta^j_{i_1 i_2
\cdots i_{l+1}}$, using the fact that if both points $x^j_{i_1 i_2 \cdots
i_{l+2}}$ exist then we have the following situation in the quiver

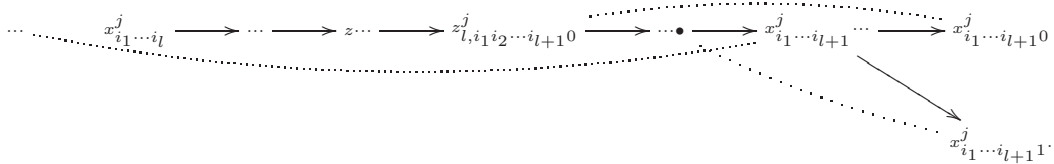
\[\xymatrix { \scriptscriptstyle \cdots  \ar@/_3ex/@{.}[rrrrrr] & \scriptscriptstyle x^j_{i_1 \cdots i_l} \ar[r] &
\scriptscriptstyle \cdots  \ar[r] & \scriptscriptstyle z \cdots
\ar[r] & \scriptscriptstyle z^j_{l,i_1 i_2 \cdots i_{l+1}0} \ar[r]
 \ar@/^2ex/@{.}[rrr]&
 \scriptscriptstyle \cdots \bullet \ar@/_1ex/@{.}[rdr] \ar[r] & \scriptscriptstyle {x^j_{i_1 \cdots i_{l+1}}}
\cdots \ar[r] \ar[dr]& \scriptscriptstyle x^j_{i_1 \cdots i_{l+1}0} \\
& & & & & & &  \scriptscriptstyle  x^j_{i_1 \cdots i_{l+1}1}. 
}
\]
The equivalence of (a) and (c) follows from the dual argument.
$\qed$

\end{Demo}

\begin{obser}
One can easily compute the top of $_A\Ext^{l+2}_A(I(c),A)$ using Corollary \ref{exists}.
\end{obser}

\section{The quiver of the higher relation extension}

{Knowing how to compute $\top\, _A \Ext_A^{i} (I(c), A)$ and
$\top \Ext_A^{i} (DA, P(z))_A$ allows us to find the new
 arrows of the higher relation extension of a string tree algebra
$A$ since they are in bijection with a basis of }

\begin{eqnarray*}
\top \Ext_A^i(DA,A)&=& \Ext_A^i(DA,A)/\rad (\Ext_A^i(DA,A))\\
&=& \displaystyle \frac{\Ext_A^i(DA,A)}{(\rad A)\Ext_A^i(DA,A) +
\Ext_A^i(DA,A) (\rad A)}.
\end{eqnarray*}

\bigskip
\bigskip

{ Note that $\Ext_A^i(DA,A). e_c= \Ext_A^i(I(c),A)$ is a left
$A$-module and that $e_z. \Ext_A^i(DA,A) = \Ext_A^i(DA,P(z))$ is a right
$A$-module.}

{Given a right (left) $A$-module $M$, we denote by $P_i(M)$ the
$i$-th term in a minimal projective resolution of $M$ and by $I_i(M)$ the
$i$-th term in a minimal injective coresolution of $M$.}

 {If we
represent the elements of $\Ext_A^i(DA,A)$ as classes $[f_{cz}]$ of
morphisms $f_{cz}: P_i(I(c)) \rightarrow P(z)$  such that the
composition of $f_{cz}$ with the map $P_{i+1}(I(c)) \rightarrow
P_i(I(c))$ of the projective resolution is zero, then we are
considering the left $A$-module structure of $\Ext_A^i(DA,A)$.
Therefore, $[f_{cz}]$ lies in $(\rad A) \Ext_A^i(DA,A)$ if and
only if  $[f_{cz}] \in \Ext_A^i(I(c),A)$ lies in the radical of
the left $A$-module $\Ext_A^i(I(c),A)$.}

{In terms of morphisms, $[f_{cz}]$ is in $\rad
\Ext_A^i(I(c),A)$ if and only if $f_{cz}$ factors non-trivially
through another morphism $f_{cy}: P_i(I(c)) \rightarrow P(y)$
such that the following diagram is commutative }

\[\xymatrix{  P_i(I(c)) \ar[rr]^{f_{cz}} \ar[dr]^{f_{cy}} & & P(z) \\
& P(y)\ar[ur]^h &  } \]

{\noindent where the map $h$ is given by the left-multiplication 
by a path in $Q$ from $z$ to $y$, and the composition of
$f_{cy}$ with $P_{i+1}(I(c)) \rightarrow P_i(I(c))$ is zero.}

\medskip

{Dually, we can represent the elements of $\Ext_A^i(DA,A)$ as
classes $[g_{cz}]$ of morphisms $g_{cz}: I(c) \rightarrow
I_i(P(z))$
 such that the composition of $g_{cz}$ with the map
$I_{i}(I(c)) \rightarrow I_{i+1}(I(c))$ of the injective
coresolution is zero. This corresponds to the right $A$-module
structure of $\Ext_A^i(DA,A)$. Therefore, $[g_{cz}]$ lies in
$\Ext_A^i(DA,A)(\rad A)$ if and only if $[g_{cz}] \in
\Ext_A^i(DA,P(z))$ lies in the radical of the left $A$-module
$\Ext_A^i(DA,P(z))$.}

{ In terms of morphisms,  $[g_{cz}]$ is in $\rad
\Ext_A^i(DA,P(z))$ if and only if $g_{cz}$ factors non-trivially
through another morphism $g_{bz}: I(b) \rightarrow I_i(P(z))$
such that the following diagram is commutative }

\[\xymatrix{  I(c) \ar[rr]^{g_{cz}} \ar[dr]^{h'} & & I_i(P(z)) \\
& I(b)\ar[ur]^{g_{bz}} &  } \]

{\noindent where the map $h'$ is given by the right-multiplication 
 by a path in $Q$ from $b$ to $c$, and the composition of
$g_{bz}$ with $I_{i}(P(c)) \rightarrow P_{i+1}(P(c))$ is zero.}

\medskip

{Moreover there is an isomorphism of vector spaces}

$$ LR: \bigoplus_{c \in Q_0} \Ext_A^i(I(c),A) \longrightarrow \bigoplus_{z \in Q_0} \Ext_A^i(DA,P(z))$$

{ \noindent such that $LR([f_{cz}])$ and $f_{cz}$ induce
the same class in $\Ext_A^i(DA,A)$. Thus, $[f_{cz}]$ is in $\rad \Ext_A^i(DA,A)$ if
and only if $[f_{cz}] \in \rad \Ext_A^i(I(c),A)$ or
$LR([f_{cz}]) \in \rad \Ext_A^i(DA,P(z))$.}

\begin{algo}\label{algo}
\ \begin{itemize}
\item[ $ \bullet$]  Compute $\top  _A\Ext^i(I(c),A)$ for all $c \in Q_0$
using Theorem \ref{left top} and Corollary \ref{exists}. (For efficiency we can restrict to the
points that are the source or the target of a relation because of Corollary 4.6 and Corollary 3.7).

\item [$ \bullet$ ] For each $c,z \in Q_0$, let $\{ \rho_{cz1}, \rho_{cz2}, \cdots \}$ be a basis for $e_z .
\top _A\Ext^i(I(c),A). $

\item [$ \bullet$] Let $B_0^i= \{ \rho_{czj}: \, \,
c,z \in Q_0 , \,c$ the source or target of relations\} be the set
that spans the vector space $\top \Ext^i(DA,A)$.

\item[$ \bullet$] Compute $\top \Ext^i(DA,P(z))_A$ for each
$z$ such that $\rho_{czj} \in B_0^i$ using Theorem \ref{left top} and the dual statements of Corollary \ref{exists}.

\item [$ \bullet$] A basis of $\top \Ext^i(DA,A)$ is

$$ B^i=B^i_0 \setminus \{\rho_{czj} \in \rad \Ext^i(DA,P(z))_A ; \, \,
c,z,j \}.$$

\end{itemize}
\end{algo}

Each element of $B^i$ has a triple subcript $czj$, and each such
element gives rise to exactly one new   arrow $z
\rightarrow c$ in the quiver of the higher relation extension.

\begin{teo} {Let $A=kQ/I$ be a string tree algebra. Then the algorithm \ref{algo} computes  two sequences
$(c_l),(z_l)$ of vertices of $Q_A$  such that the arrows in the quiver of the higher relation extension  are exactly those of $Q_A$ plus one
additional arrow from each $z_l$ to $c_l$. }
\end{teo}
\begin{Demo}
This follows from the discussion preceding the algorithm.\qed
\end{Demo}

\begin{obser}
The vertices $(c_l),(z_l)$
are not necessarily distinct, there may be repetitions.
\end{obser}

\begin{ejem} Let $A=kQ/I$ be the string algebra   given by the following bound quiver:

\[\xymatrix@R12pt@C12pt{\scriptscriptstyle1\ar@/^2ex/@{.}[rrr]\ar[r]&\scriptscriptstyle2\ar[r]&
\scriptscriptstyle3\ar[r]&\scriptscriptstyle4.}
\]
Then there exists an element $\rho_{2,4}\in e_4. \top _A\Ext^2_A(I(2),A)$ which is not in $\top \Ext^2_A(DA,P(4))_A . e_2$ and therefore not in $ \top_A\Ext^2_A(DA,A)_A$.
Thus the quiver of the higher relation extension is 
\[\xymatrix@R12pt@C12pt{\scriptscriptstyle1\ar[r]&\scriptscriptstyle2\ar[r]&
\scriptscriptstyle3\ar[r]&\scriptscriptstyle4.\ar@/_2ex/[lll]}
\]

\end{ejem}

\begin{ejem} In this example, the higher relation extension contains an $\Ext^2$-arrow $5\to 1$ although there is no relation between the points 5 and 1.
Let $A=kQ/I$ be the string algebra   given by the bound quiver:

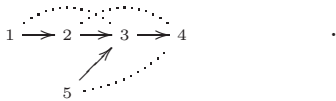
\[\xymatrix@R12pt@C12pt{  \scriptscriptstyle 1  \ar[r] \ar@/^2ex/@{.}[rr] &  \scriptscriptstyle 2 \ar[r] \ar@/^2ex/@{.}[rr]&
 \scriptscriptstyle 3 \ar[r] &  \scriptscriptstyle 4 & & & . \\
&  \scriptscriptstyle 5 \ar[ru] \ar@/_1ex/@{.}[rur] } \]

\noindent Then the quiver of the higher relation extension of $A$ is
the following:

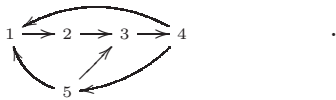
\[\xymatrix@R12pt@C12pt{  \scriptscriptstyle 1  \ar[r]  &  \scriptscriptstyle 2 \ar[r] &
 \scriptscriptstyle 3 \ar[r] &  \scriptscriptstyle 4 \ar@/_2ex/[lll] \ar@/^1ex/[dll] & & & .\\ &  \scriptscriptstyle 5 \ar[ru] \ar@/^1ex/[ul]} \]

\end{ejem}

\begin{ejem} Let $A$ be the string algebra of Example \ref{example 1}.
Then the quiver of the higher relation extension of $A$ is the following:

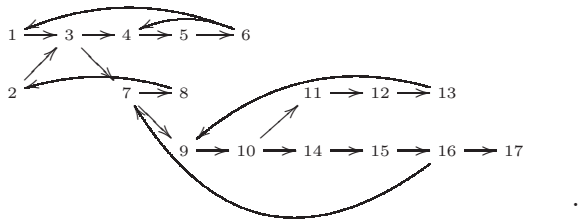
\[\xymatrix@R12pt@C12pt{ \scriptscriptstyle 1 \ar[r] & \scriptscriptstyle 3  \ar[dr] \ar[r] &
 \scriptscriptstyle 4 \ar[r] &\scriptscriptstyle 5  \ar[r] & \scriptscriptstyle 6 \ar@/_2ex/[llll] \ar@/_/[ll]\\
\scriptscriptstyle 2 \ar[ur] & & \scriptscriptstyle 7 \ar[r]
\ar[dr]&  \scriptscriptstyle 8 \ar@/_/[lll]& & \scriptscriptstyle 11
\ar[r]  &
 \scriptscriptstyle 12 \ar[r] & \scriptscriptstyle 13 \ar@/_3ex/[lllld]\\
& & &  \scriptscriptstyle 9
 \ar[r] & \scriptscriptstyle 10 \ar[r] \ar[ur] & \scriptscriptstyle 14 \ar[r] &
 \scriptscriptstyle 15 \ar[r] & \scriptscriptstyle 16 \ar@/^7ex/[dllllluu] \ar[r] & \scriptscriptstyle 17 \\
 & & & & & & & & &.}
\]
\end{ejem}

\begin{ejem} This example illustrates the situation in Corollary 4.6 (b). Let $A=kQ/I$ be the string algebra   given by the bound quiver:

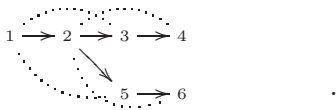
\[\xymatrix@R12pt@C12pt{\scriptscriptstyle 1\ar@/^2ex/@{.}[rr] \ar@/_2ex/@{.}[drr]\ar[r] &\scriptscriptstyle 2\ar@/^2ex/@{.}[rr] \ar@/_3ex/@{.}[drr]\ar[r] \ar[dr]&
\scriptscriptstyle 3\ar[r]&\scriptscriptstyle 4 \\
&  & \scriptscriptstyle 5\ar[r]&\scriptscriptstyle 6 & & & .} \]

Then the quiver of the higher relation extension of $A$ is the following:

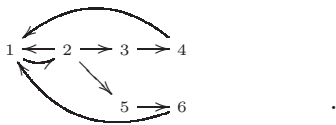
\[\xymatrix@R12pt@C12pt{\scriptscriptstyle 1 \ar@/_1ex/[r] &\scriptscriptstyle 2 \ar[l]\ar[r] \ar[dr]&
\scriptscriptstyle 3\ar[r]&\scriptscriptstyle 4 \ar@/_3ex/[lll]\\
&  & \scriptscriptstyle 5\ar[r]&\scriptscriptstyle 6 \ar@/^3ex/[ulll]& & & .} \]
Note the existence of a 2-cycle.
\end{ejem}

\section{The higher relation bimodule for gentle algebras}

{ Recall that  a set of monomial relations
$\{\kappa_i\}_{i=1,..,t}$ is called an {\bf overlapping} if the
paths $\kappa_i$ and $\kappa_{i+1}$ have a common subpath  $
\vartheta$ such that $\kappa_i= \vartheta_i  \vartheta$ and
$\kappa_{i+1}= \vartheta \vartheta_{i+1}$, for all $i=1,..,t-1$. A
 {\bf maximal $t$-overlapping } is an overlapping $\{\kappa_i\}_{i=1,..,t}$
such that there exists no monomial relation $\kappa$ such that the
sets $\{ \kappa, \kappa_i, \, \, i=1, \cdots, t \}$ and  $\{
\kappa_i, \, \, i=1, \cdots, t, \kappa \}$ are an overlapping, see \cite{GHZ,Gu}.}

\begin{lema} Let $\kappa= (\kappa_1, \cdots, \kappa_t)$ be the following maximal
$t$-overlapping over a gentle algebra $A=kQ/I$:

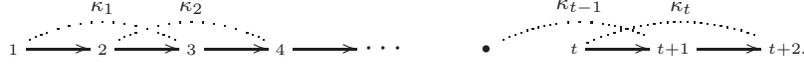
\[ \xymatrix{  \scriptscriptstyle 1  \ar[r] \ar@/^2ex/@{.}[rr]^{\kappa_1} &  \scriptscriptstyle 2 \ar[r] \ar@/^2ex/@{.}[rr]^{\kappa_2}&
 \scriptscriptstyle 3 \ar[r] &  \scriptscriptstyle 4  \ar[r] & \cdots &\scriptscriptstyle \bullet  \ar@/^2ex/@{.}[rr]^{\kappa_{t-1}}
 & \scriptscriptstyle t \ar[r] \ar@/^2ex/@{.}[rr]^{\kappa_t}
 &  \scriptscriptstyle t+1 \ar[r]   & \scriptscriptstyle t+2 .} \]

\noindent Then, for the injective $I(1)$ associated to the vertex
$1$, the sequence of $x_{i_1i_2\cdots i_t}$ is:

$$x_0=3, \, x_{00}=4,  \, x_{000}=5,  \, \cdots, x_{i_1i_2\cdots i_t}=
x_{00 \cdots 0}= t+2.$$

\end{lema}

\begin{Demo} This follows from the construction of the points $x_{i_1i_2\cdots
i_t}$ given in section 3. $\qed$

\end{Demo}

\begin{obser} Observe that there may be other points $x_{i_1i_2\cdots
i_t}$ where some $i_j \ne 0$. In the Lemma we only consider  
one branch of the quiver which contains all the points $x_{00\cdots 0}$.
\end{obser}

\begin{prop}\label{prop 6.3} For every maximal $t$-overlapping $\kappa= (\kappa_1, \cdots, \kappa_t)$ from $c$ to $z$ there is exactly
one new arrow $\alpha(\kappa): z \rightarrow c$ in the higher relation extension
which is induced by an element of $\Ext_A^{l+1}(I(c),
P(z))$ and these are the only new arrows in the higher relation extension.
Moreover, we have the following relations:
\begin{itemize}
\item[(a)] $\alpha(\kappa) \alpha_1=0$  and $\alpha_{t+1}\alpha_{\kappa}=0$, where $\alpha_1$ and $\alpha_{t+1}$ denote the first and
the last arrow of $\kappa$;
\item[(b)] $\zeta \rho \zeta'$ where $\zeta,
\zeta'$ are new arrows and $\rho$ is a path consisting of old arrows.
\end{itemize}
\end{prop}

\begin{Demo} By Corollary 4.5, $\Ext_A^{l+2}(I(c), P(z)) \ne 0$ if and only
if there is a non-zero path $\omega: z \rightsquigarrow
x_{i_1i_2\cdots i_{l+1}} $ not passing through $x_{i_1i_2\cdots
i_{l}}$ and such that the compositions with the non-zero paths
$\omega'_{i_{l+2}}: x_{i_1i_2\cdots i_{l+1}} \rightsquigarrow
x_{i_1i_2\cdots i_{l+2}}$ are both zero if $i_{l+2}$ exists, see
figure.

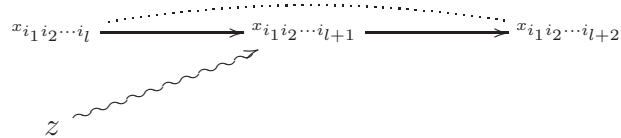
\[\xymatrix{   & \scriptscriptstyle x_{i_1i_2\cdots i_{l}}  \ar[rr] \ar@/^2ex/@{.}[rrrr]
& & \scriptscriptstyle x_{i_1i_2\cdots i_{l+1}} \ar[rr]  & &
 \scriptscriptstyle x_{i_1i_2\cdots i_{l+2}} \\
&  z \ar@{~>}[urr]& && & & &.} \]

But the previous Lemma implies that 
$x_{i_1i_2\cdots i_{l}} \rightarrow x_{i_1i_2\cdots i_{l+1}}
\rightarrow x_{i_1i_2\cdots i_{l+2}}$ is a relation of length 2, contradicting that $A$ is gentle. Therefore $i_{l+2}$ does not exist, that is, $\pd I(c)=l+2$.
Then we have the situation

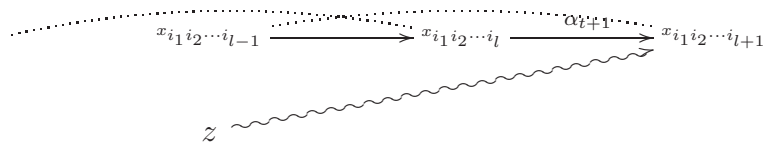
\[\xymatrix{  & \ar@/^2ex/@{.}[rrrr] & & \scriptscriptstyle x_{i_1i_2\cdots i_{l-1}}  \ar[rr] \ar@/^2ex/@{.}[rrrr]
& &
 \scriptscriptstyle x_{i_1i_2\cdots i_{l}} \ar[rr]^{\alpha_{t+1}} &&
 \scriptscriptstyle x_{i_1i_2\cdots i_{l+1}} \\
& & & z \ar@{~>}[urrrr]} \]

\noindent and $x_{i_1i_2\cdots i_{l+1}}$ is the target of an
overlapping $\omega$.

\noindent Thus by Corollary 4.5,
$\Ext_A^{l+2}(I(c),P(z)) \ne 0$ if and only if there is a non-zero
path $\omega$ from $z$ to $x_{i_1i_2\cdots i_{l+1}}$ not passing
through $x_{i_1i_2\cdots i_{l}}$. Then, by Corollary
\ref{exists} a2), $\omega$ induces an element of the top of
$_A\Ext_A^i(I(c),A)$ if and only if $z=x_{i_1i_2\cdots i_{l+1}}$.

\noindent To check whether $\omega: z \rightarrow x_{i_1i_2\cdots
i_{l+1}}$ induces an element of $\top \Ext_A^{l+2}(DA,A)$, we
can apply the algorithm \ref{algo}. Hence, $\omega$ induces an
element of the $\top \Ext_A^{l+2}(DA,A)$ if and only if $z$ is the
starting point of the overlapping $\kappa$.

\noindent The result about the new arrows now follows from the
algorithm.

Using the fact that $\Ext_A^{l+2}(I(c),P(z)) \ne 0$ if and only if there is a non-zero
path $\omega$ from $z$ to $x_{i_1i_2\cdots i_{l+1}}$ not passing
through $x_{i_1i_2\cdots i_{l}}$ with $z=x_{i_1i_2\cdots i_{l}}$ shows that  
 $\alpha_{t+1}
\alpha(\kappa)=0$. Dually, one proves that $\alpha(\kappa)\alpha_1=0$, and the relations of the form $\zeta\rho\zeta'$ occur since we are dealing with a trivial extension.
 $\qed$
\end{Demo}

The following example shows that the higher relation extension of a gentle algebra is not necessarily gentle.

\begin{ejem}
Let $A$ be given by the bound quiver
\[\xymatrix{ 1\ar[r] \ar@{.}@/^10pt/[rr]&2\ar[r]&3&4\ar[l]^\rho\ar[r] \ar@{.}@/^10pt/[rr] &5 \ar[r]&6.
}
\]
Then the higher relation extension coincides with the relation extension and has the quiver 
\[\xymatrix{ 1\ar[r] &2\ar[r]&3\ar@/_10pt/[ll]_{\zeta'}&4\ar[l]_\rho\ar[r]  &5 \ar[r]&6\ar@/_10pt/[ll]_{\zeta}
}
\]
bound by relations of length 2 and the relation $\zeta\rho\zeta'$, which is of length 3.
\end{ejem}

\begin{coro}
\label{tensor rel}
The tensor algebra of the higher relation bimodule has the same quiver as the higher relation extension and has the relations in Proposition \ref{prop 6.3} (a). In particular its relation ideal is quadratic. \qed
\end{coro}

\section{The tensor algebra of a gentle algebra}

\begin{teo} Let $A$ be a gentle algebra.   
\begin{itemize}
\item[(a)] The tensor algebra
$T_A(\bigoplus_{i \geq 2} \Ext_A^i(DA,A))$ is gentle.
\item[(b)] The higher relation extension $A \ltimes (\bigoplus_{i \geq 2} \Ext_A^i(DA,A))$ is monomial. 
\end{itemize}
\end{teo}

\begin{Demo} Since the universal cover of a gentle algebra is a
gentle tree, we may assume that $A$ is a tree. We prove the  
conditions S1), S2), S3), G1) and G2) of section 2.

\begin{itemize}
\item[S2)] At every point there are at most two incoming
arrows (dually, outcoming arrows).

Suppose there are three arrows $\alpha, \beta, \gamma$ with target $x$ (see figure)

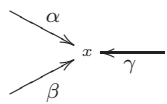
\[\xymatrix@R10pt {  \ar[rd]^{\alpha}&  & \\
 & \scriptscriptstyle x  & \ar[l]^{\gamma} \\
\ar[ru]_{\beta}  & & } \]

 Then at least one, say $\gamma$, is a new arrow. Hence, $\gamma$ corresponds to an overlapping
$\omega=(\omega_1, \omega_2, \cdots)$ with source $x$ and there
is no relation involving $\alpha$ or $\beta$ and overlapping with
$\omega_1$, that is,
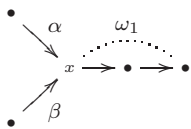
\[\xymatrix@R12pt@C12pt{   \scriptscriptstyle \bullet \ar[rd]^{\alpha}&  & & \\
 & \scriptscriptstyle x  \ar[r] \ar@/^2ex/@{.}[rr]^{\omega_1} &  \scriptscriptstyle \bullet \ar[r]
 & \scriptscriptstyle \bullet  \\
\scriptscriptstyle \bullet \ar[ur]_{\beta}  & & &} \]

Because $A$ is gentle, at least one of the arrows $\alpha$ and $\beta$ is new. Assume $\alpha $ is old and $\beta$ is new.
 Then we have
two such overlappings $\omega, \, \omega'$ and no relation
involving $\alpha$ and overlapping with $\omega$ or $\omega'$, that
is we have the following situation in the bound quiver of $A$.
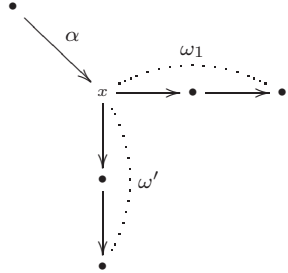
\[\xymatrix {   \scriptscriptstyle \bullet \ar[rd]^{\alpha}&  & & \\
 & \scriptscriptstyle x  \ar[r] \ar[d] \ar@/^2ex/@{.}[rr]^{\omega_1} \ar@/^2ex/@{.}[dd]^{\omega'} &  \scriptscriptstyle \bullet \ar[r]
 & \scriptscriptstyle \bullet  \\
 & \scriptscriptstyle \bullet \ar[d] & & \\
  & \scriptscriptstyle \bullet & &} \]
which yields a contradiction.

Finally, if all   three arrows $\alpha, \beta, \gamma$ are new, we
get three overlappings starting at $x$. Because $A$ is gentle, condition G1) implies that we have three arrows having $x$ as a source, a contradiction.

\item[S1,G2)] Suppose we have a minimal relation involving at least two
paths in the sense of \cite{MP}. Then, in the higher relation extension we have
at least two paths $c_1, c_2$ starting and ending at the same point
with at least one new arrow in each of these paths. Let $c_i= c_{i1}
\alpha_i c_{i2}$ where $\alpha_i$ is a new arrow, $i=1,2$.

Assume first that there is exactly one new arrow on each path $c_i$. Then
each $\alpha_i$ corresponds to an overlapping $\omega_i$ in $A$
starting at the target of $\alpha_i$ and ending at its source, and this
contradicts the assumption that $A$ is a tree.

If $c_i$ contains several new arrows, the same argument as before
applies.

This shows that  the higher relation extension of  $A$ is monomial and hence that the tensor algebra is also monomial, and even has a quadratic relation ideal, because of Corollary \ref{tensor rel}.

\item[S3)] Suppose we have the following subquiver
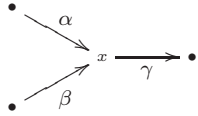
\[\xymatrix@R10pt{ &\scriptscriptstyle \bullet  \ar[dr]^\alpha\\
& & \scriptscriptstyle x  \ar[r]_\gamma & \scriptscriptstyle \bullet\\
& \scriptscriptstyle \bullet \ar[ur]_\beta &} \] such that $\alpha
\gamma$ and $\beta \gamma$ are not in the ideal $I$ of the tensor
algebra $T_A(\bigoplus_{i \geq 2} \Ext_A^i(DA,A))$. Then one of the
three arrows is new. First assume that $\gamma$ is a new arrow, then
$\gamma$ corresponds to an overlapping $\omega$ ending at $x$ which
implies that $\alpha$ or $\beta$ must be a new arrow, say $\beta$
which correspond to another overlapping $\omega'$ starting at $x$.
But then the last arrow of $\omega$ and the first arrow of $\omega'$
are not bound by a relation and also $\alpha$ is not bound by a relation
with the first arrow of $\omega'$. This contradicts $A$ being
gentle.

Suppose now that $\alpha$ is a new arrow corresponding to an
overlapping $\omega$ starting at $x$. Because of the first case, we may
assume that $\gamma$ is not new. Since $\beta$ is not bound by a relation with
the first arrow in $\omega$, it must be with $\gamma$,
contradicting the assumption $\beta \gamma \notin I$.

\item [G1)] Suppose we have a subquiver
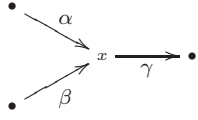
\[\xymatrix@R10pt{ &\scriptscriptstyle \bullet  \ar[dr]^\alpha\\
& & \scriptscriptstyle x  \ar[r]_\gamma & \scriptscriptstyle \bullet\\
& \scriptscriptstyle \bullet \ar[ur]_\beta &} \] such that $\alpha
\gamma$ and $\beta \gamma$ are in the relation ideal of the tensor
algebra $T_A(\bigoplus_{i \geq 2} \Ext_A^i(DA,A))$. If $\gamma$ is a
new arrow corresponding to an overlapping $\omega_{\gamma}$ ending
at $x$ then $\alpha$ or $\beta$ must be new, say $\beta$,
and corresponding to an overlapping $\omega_{\beta}$ as above, which is bound by no relation with $\alpha$. It follows from our description of the
bound quiver that the new arrow $\beta$ is not bound by a relation
with $\gamma$, because $\gamma$ is not in the overlapping
$\omega_{\beta}$, and this is a contradiction. $\qed$

\end{itemize}
\end{Demo}

ACKNOWLEDGMENTS:
The first author gratefully acknowledges support from the NSERC of Canada, the FQRNT of Qu\'ebec and the Universit\'e de Sherbrooke, the second author gratefully acknowledges support from the NSERC of Canada and the CONICET of Argentina, and the third author gratefully acknowledges support from NSF grant DMS-1001637 and the University of Connecticut. 


\small
\noindent I. Assem, D\'epartement de math\'ematiques, Universit\'{e} de Sherbrooke, Sherbrooke, Qu\'{e}bec, Canada, JIK2R1. \\
Ibrahim.Assem@usherbrooke.ca\\
M.A. Gatica, Departamento de Matem\'atica, Universidad Nacional del Sur, Avenida Alem 1253, (8000) Bah\'{\i}a Blanca, Buenos Aires, Argentina. \\
mariaandrea.gatica@gmail.com\\
R. Schiffler, Department of Mathematics, 196 Auditorium Road, Storrs, CT 06269-3009, University of Connecticut, Connecticut, USA.\\
schiffler@math.uconn.edu

\end{document}